\documentclass[12pt]{amsart}
\usepackage{
  a4wide, 
  amsmath,
  amsfonts,
  amssymb,
  graphicx, 
  times,
  color,
  hyphenat,
  pinlabel,
  stmaryrd, mathabx, scalefnt, rotating 
}

\usepackage[french]{babel}
\usepackage[utf8]{inputenc}
\usepackage[T1]{fontenc}
\addto{\captionsfrench}{}
\input xy
\xyoption{all}
\DeclareMathOperator{\fun}{Fun}
\DeclareMathOperator{\inv}{Inv}
\DeclareMathOperator{\rot}{rot}

\newcommand{\figins}[3] 
{\raisebox{#1pt}{\includegraphics[height=#2 in]{figs/#3}}}

\newtheorem{thm}{Theorem}[section]

\theoremstyle{definition}
\newtheorem{defn}[thm]{Definition}

\DeclareMathOperator{\Ind}{Ind}
\DeclareMathOperator{\Res}{Res}

\catcode`\@=11
\long\def\@makecaption#1#2{%
    \vskip 10pt
    \setbox\@tempboxa\hbox{%
\small{#1: }\ignorespaces #2}%
    \ifdim \wd\@tempboxa >\captionwidth {%
        \rightskip=\@captionmargin\leftskip=\@captionmargin
        \unhbox\@tempboxa\par}%
      \else
        \hbox to\hsize{\hfil\box\@tempboxa\hfil}%
    \fi}
\newdimen\@captionmargin\@captionmargin=2\parindent
\newdimen\captionwidth\captionwidth=\hsize
\catcode`\@=12
\title{On Jaeger's HOMFLY-PT expansions, branching rules and link homology: a progress report
}
\author{Pedro Vaz}
\address{Institut de Recherche en Math\'ematique et Physique\\
Universit\'e Catholique de Louvain\\ 
Chemin du Cyclotron 2\\ 
B 1348 Louvain-la-Neuve\\ 
Belgium}
\email{pedro.vaz@uclouvain.be}
\keywords{Quantum invariant, branching rules, link homology.}
\begin{document}
%
%
\newdimen\captionwidth\captionwidth=\hsize
%
%
\begin{abstract}
  We describe Jaeger's HOMFLY-PT expansion of the Kauffman polynomial and how to generalize it  
to other quantum invariants using the so-called ``branching rules'' 
for Lie algebra representations. We present a program which aims to construct  
Jaeger expansions for link homo\-lo\-gy theories.
This note is an updated write-up of a talk given by the author at the 
Meeting of the Sociedade Portuguesa de Matem\'atica in July 2012. 
\end{abstract}
\maketitle

%
%
\pagestyle{myheadings}
\markright{Pedro Vaz}
\markboth{\em\small Pedro Vaz}{\em\small Jaeger's HOMFLY-PT expansions, branching rules and link homology}
%
%
%
\section{Link polynomials and Jaeger expansions}
\label{sec:polys}

This story starts with two celebrated invariant polynomials of links.
\begin{defn}
The \emph{Kauffman polynomial} $F=F(a,q)$ is the unique invariant of framed uno\-riented links 
satisfying
\begin{gather*}
F\left(
\figins{-8}{0.3}{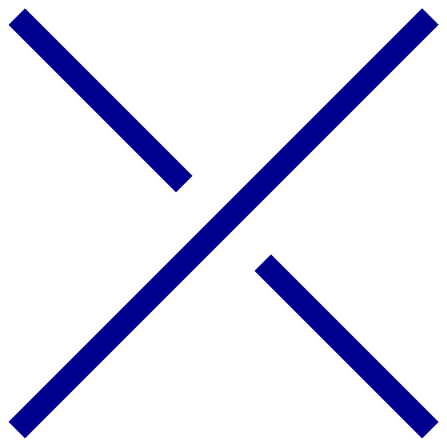}\right) \  -  \ 
F \left(
\figins{-8}{0.3}{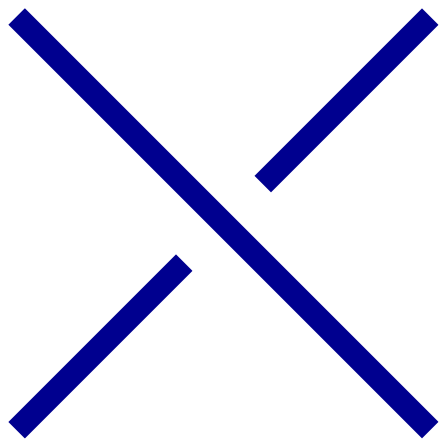}\right)\
=\  (q - q^{-1})
\biggl(\ 
F\left(
\figins{-8}{0.3}{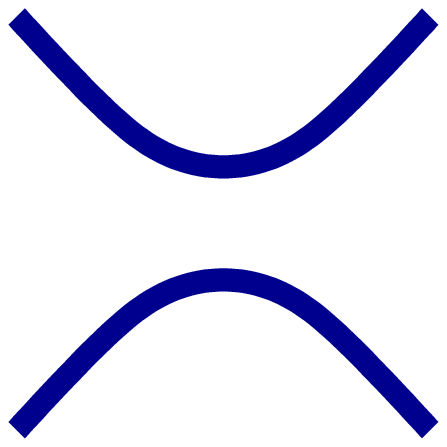}\right)\  -  \
F\left(
\figins{-8}{0.3}{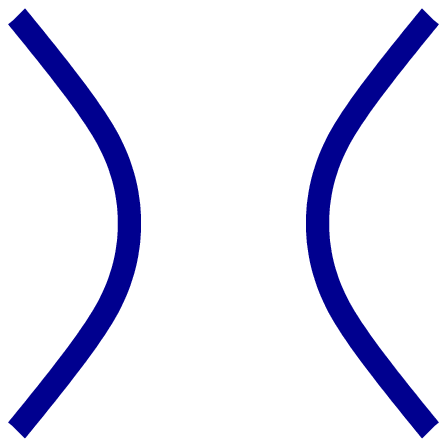}\right)\
\biggr)
\\[1ex]
F\left(\;
\figins{-9}{0.35}{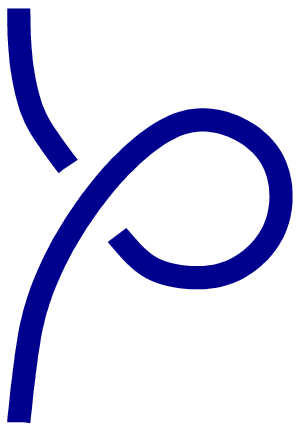}\right)\  
= \ a^{-2}q\  
F\left(\ 
\figins{-9}{0.35}{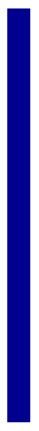}\ \right)
\qquad\text{and}\qquad
F\left(\figins{-6}{0.25}{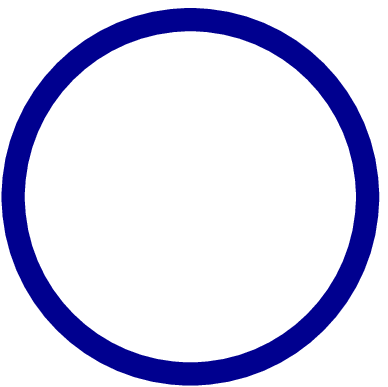}\right)\  
=\  \frac{a^2q^{-1}-a^{-2}q}{q-q^{-1}}+1 .
\end{gather*}
\end{defn}
\begin{defn}
The \emph{HOMFLY-PT polynomial} $P=P(a,q)$ is the unique invariant of framed oriented  links satisfying
\begin{gather*}
a P\left(
\figins{-8}{0.3}{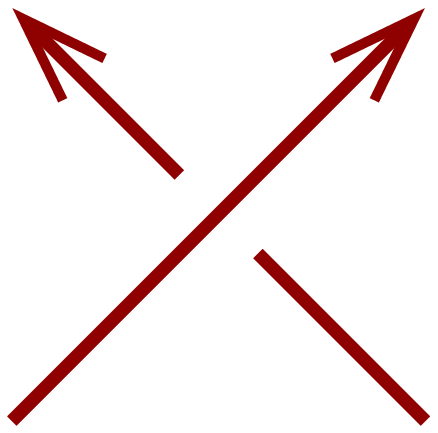}\right)  \ - \  
a^{-1}P\left(
\figins{-8}{0.3}{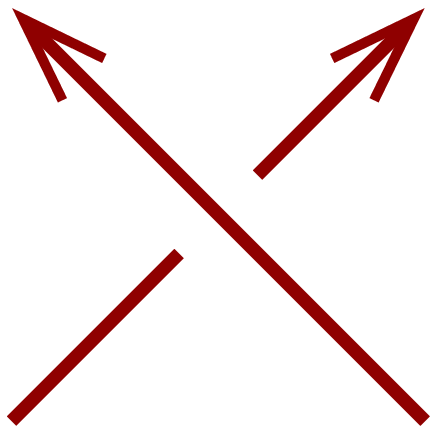}\right)
=\  (q - q^{-1})P \left(
\figins{-8}{0.3}{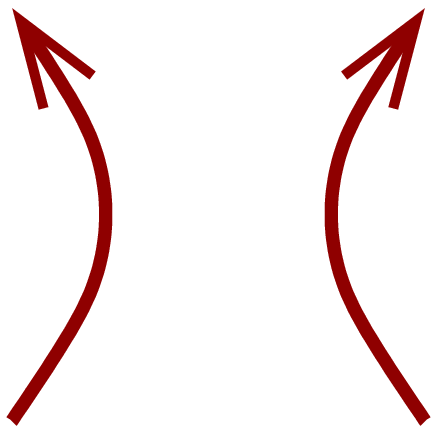}
\right)
\\[1ex]
P\left(
\figins{-9}{0.35}{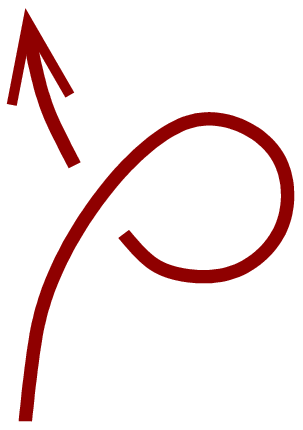}\right)
=\ a^{-1}
P\left(
\figins{-9}{0.35}{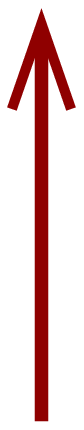}\right) , 
\mspace{30mu}
P\left(
\figins{-9}{0.35}{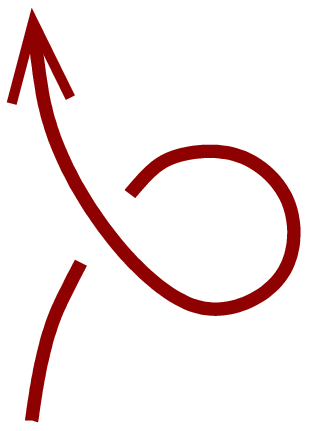}\right)
=\ a
P\left(
\figins{-9}{0.35}{one-u}\right)
\quad\text{and}\quad
P\left(\figins{-6}{0.25}{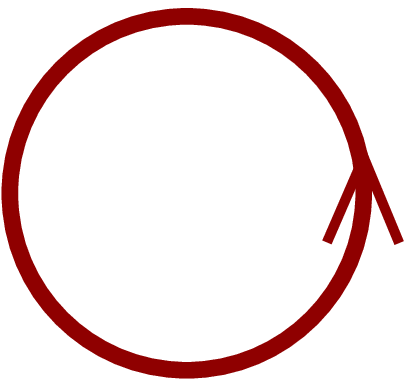}\right)\  
= \ \frac{a - a^{-1}}{q-q^{-1}}  .
\end{gather*}
\end{defn}

\medskip

Either the Kauffman and the HOMFLY-PT polynomials can be normalized to give 
ambient isotopy invariants of (oriented) links, but we will not proceed in this direction.   
In 1989 François Jaeger showed that the Kauffman polynomial of a link $L$ can be obtained as a 
weighted sum of HOMFLY-PT polynomials on certain links associated to $L$. 
Consider the following formalism 
\begin{equation}\label{eq:jxpan}
\begin{aligned}
\left[
\figins{-8}{0.30}{Xing-p}\right]    &=   
(q-q^{-1})
\biggl( \left[
\figins{-8}{0.30}{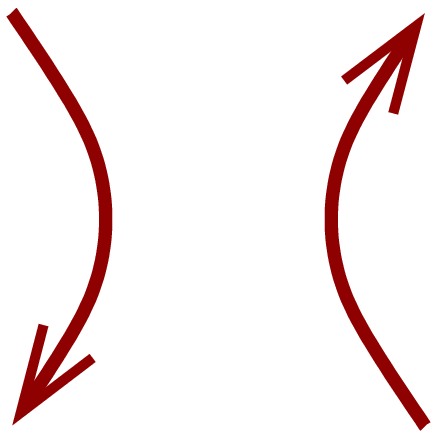}\right]    -   
\left[
\figins{-8}{0.30}{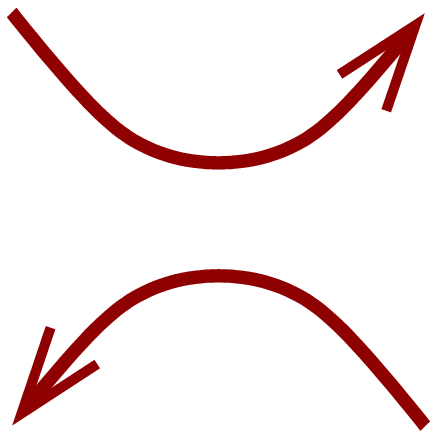}\right] \biggr)   
+   
\left[\figins{-8}{0.30}{Xing-pu}\right]   +  
\left[\figins{-8}{0.30}{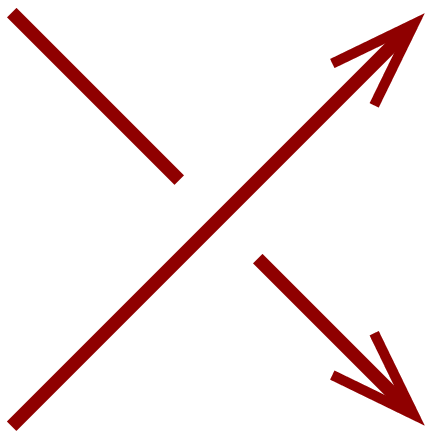}\right]   +  
\left[\figins{-8}{0.30}{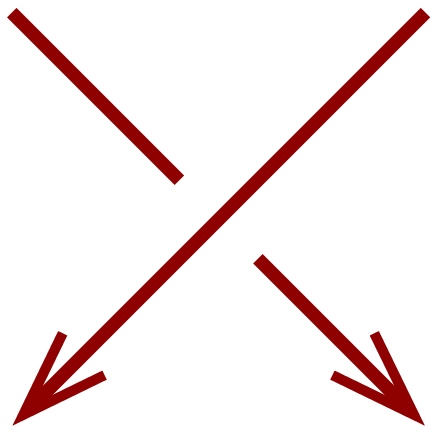}\right]   +  
\left[\figins{-8}{0.30}{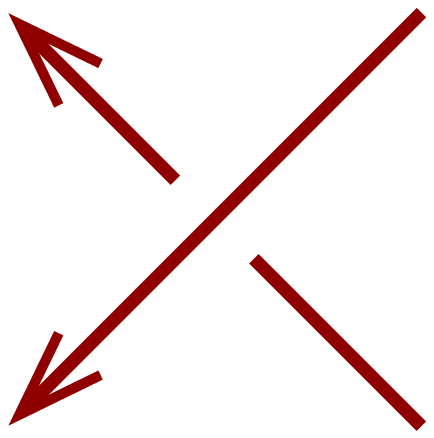}\right]
\\[1ex]
\left[\figins{-8}{0.30}{loop}\right]\    &=   
\left[\figins{-8}{0.30}{loop-r}\right]  +   
\left[\figins{-8}{0.30}{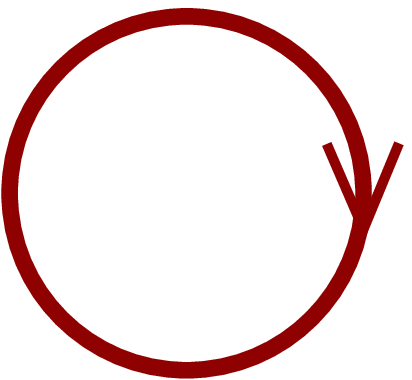}\right]
\end{aligned}
\end{equation}
where the r.h.s. is evaluated to HOMFLY-PT polynomials completed with  
information about \emph{rotation numbers}:
$[\vec{D}] = (a^{-1}q)^{\rot D} P(\vec{D})$ 
(of course we only take the diagrams which are globally coherently oriented).
The expression in~\eqref{eq:jxpan} can be seen as a state expansion with the coefficients given as vertex weights. 
For a state $\sigma$ we denote its weight by $w(\sigma)$. 
The proof of the following can be found in~\cite{kauffman}.
\begin{thm}[F.~Jaeger, 1989]\label{thm:jaeg}
Let $D$ be a diagram of a link $L$ and $\Sigma(D)$ denote the set of states of $D$. 
The sum
\begin{equation}\label{eq:jXpan}
\sum\limits_{\sigma\ \in\ \Sigma(D)}w(\sigma)[\sigma] = F(D)
\end{equation}
is a HOMFLY-PT expansion of the Kauffman polynomial of $L$.
\end{thm}

It is not hard to see how to extend this expansion to tangles (see~\cite{VW}). 
To explain this expansion we look into the representation theory of quantum enveloping algebras 
of the simple Lie algebras (QEAs). 
It is known that the HOMFLY-PT polynomial is related to the representation theory 
of the QEAs of type $A_{n-1}$ (e.g. $\mathfrak{sl}_{n}$) and that in turn 
the Kauffman polynomial is related to the QEAs of types $B_{n}$, $C_n$ and $D_n$ 
($\mathfrak{so}_{2n+1}$, $\mathfrak{sp}_{2n}$ and $\mathfrak{so}_{2n})$ respectively).
Taking $a=q^n$ in $F(a,q)$ and $P(a,q)$ we obtain the $\mathfrak{so}_{2n}$ and 
the $\mathfrak{sl}_n$ polynomials respectively.

\medskip

Following N. Reshetikhin and V. Turaev there is a functor from the category of tangles whose   
arcs are colored by irreducible finite dimensional (f.d.) representations of a  
QEA $\mathfrak{g}$ to the (tensor) category of f.d. representations
of $\mathfrak{g}$ (see~\cite{resh-tur, tur}). 
In other words, for each of these tangles there is a $\mathfrak{g}$-invariant map which depends only on the 
(regular) isotopy class of the tangle which 
gives a full isotopy invariant in the cases we are interested in. 
So what we really have in Theorem~\ref{thm:jaeg} is an $\mathfrak{sl}_n$-expansion of 
the $\mathfrak{so}_{2n}$-polynomial 
(the case where all strands are colored by the fundamental representation)!

\medskip

There are other (2-variable) HOMFLY-PT expansions of $F(a,q)$ resulting 
in $\mathfrak{sl}_n$-expansions of the  
$\mathfrak{so}_{2n+1}$ and $\mathfrak{sp}_{2n}$ polynomials after the specialization $a=q^n$.
For example the assignment (this was found by the author together with E. Wagner 
and will appear somewhere in the literature~\cite{VW2})
\begin{align*}
\left[\figins{-8}{0.3}{Xing-p}\right]_{\! B_n}\!\!\!  = &
(q-q^{-1})\! \biggl(\!\!
\left[\figins{-8}{0.3}{downup}\right]_{\! B_n}\!\!\! + 
\left[\figins{-8}{0.3}{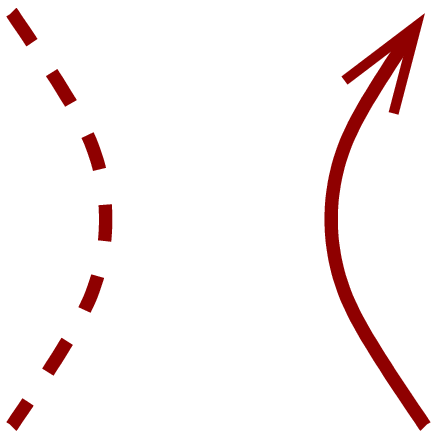}\right]_{\! B_n}\!\!\!  + 
\left[\figins{-8}{0.3}{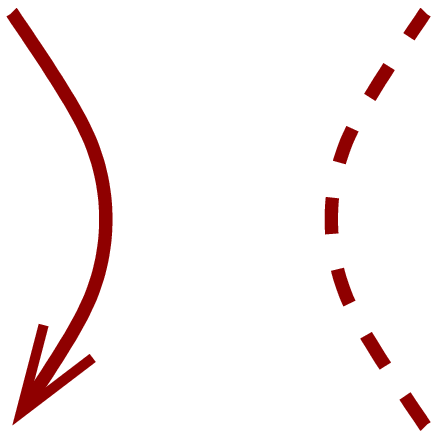}\right]_{\! B_n}\!\!\! - 
\left[\figins{-8}{0.3}{cupcap-lr}\right]_{\! B_n}\!\!\! - 
\left[\figins{-8}{0.3}{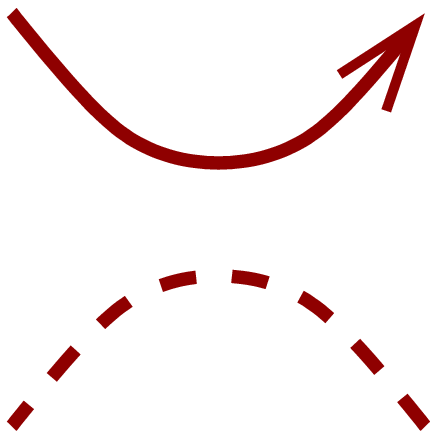}\right]_{\! B_n}\!\!\! - 
\left[\figins{-8}{0.3}{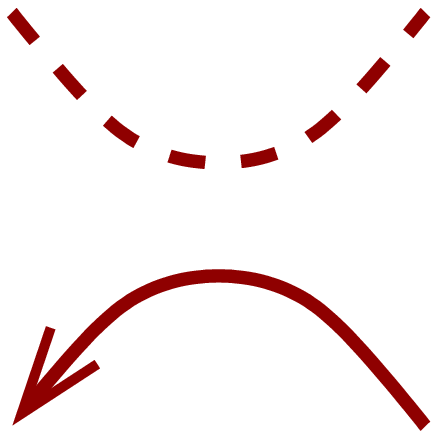}\right]_{\! B_n}\!  
\biggr)
\\[1ex] & 
  + 
\left[\figins{-8}{0.3}{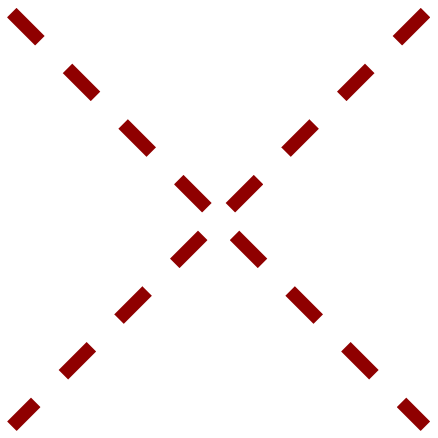}\right]_{\! B_n}\!\!\! +  
\left[\figins{-8}{0.3}{Xing-pu}\right]_{\! B_n}\!\!\! +  
\left[\figins{-8}{0.3}{Xing-pd}\right]_{\! B_n}\!\!\! +  
\left[\figins{-8}{0.3}{Xing-lp}\right]_{\! B_n}\!\!\! +  
\left[\figins{-8}{0.3}{Xing-rp}\right]_{\! B_n}\!\!\!
 -\left[\figins{-8}{0.3}{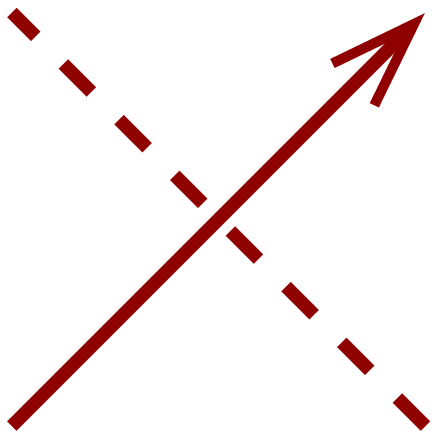}\right]_{\! B_n}\!\!\!  - 
\left[\figins{-8}{0.3}{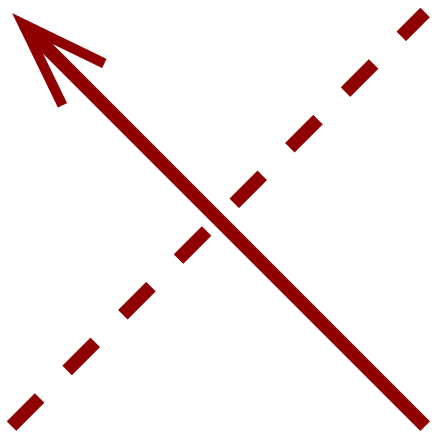}\right]_{\! B_n}\!\!\!  
\\[1ex] &
- 
\left[\figins{-8}{0.3}{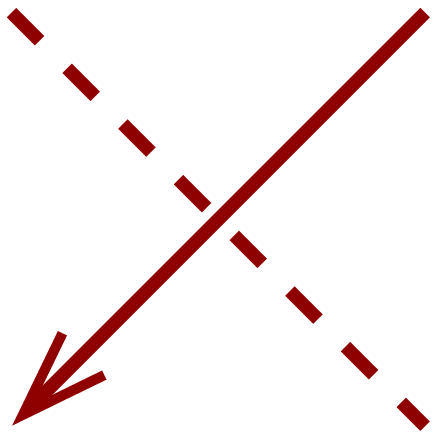}\right]_{\! B_n}\!\!\!  - 
\left[\figins{-8}{0.3}{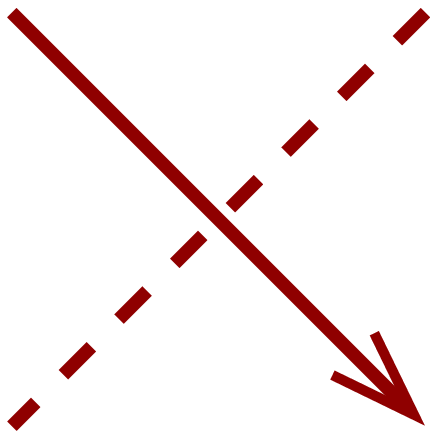}\right]_{\! B_n} 
\\[1.5ex]
\left[ \ \figins{-8}{0.3}{one}\  \right]_{B_n}  =&  
\left[ \ \figins{-8}{0.3}{one-u}\  \right]_{B_n}  +    
\left[ \ \figins{-8}{0.3}{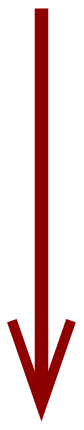}\  \right]_{B_n}   +   
\left[ \ \figins{-8}{0.3}{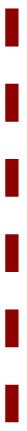} \ \right]_{B_n}
\end{align*}
where  
$[\vec{D}]_{B_n}=a^{-\rot D}P(\vec{D})$ 
and a dashed line means the corresponding strand is to be erased, 
gives an  
$\mathfrak{sl}_n$-expansion of the $\mathfrak{so}_{2n+1}$-polynomial.

\section{Branching rules, link homology and categorification}
\label{sec:categ}

Let us give an explanation for this phenomenon. An 
inclusion $\mathfrak{l}\hookrightarrow \mathfrak{g}$ of Lie algebras (resp. QEAs)
gives rise to functors ($\Ind$ and $\Res$) between their 
categories of representations. 
In general $\Res$ does not send an irreducible $M$ over $\mathfrak{g}$ to an 
irreducible over $\mathfrak{l}$ but if we restrict ouselves to finite-dimensional representations we know that 
$\Res(M)$ decomposes as a direct sum of irreducibles for $\mathfrak{l}$. 
The branching rule tells us how to obtain such a decomposition \emph{i.e.} how to express an irreducible for $\mathfrak{g}$ 
as a direct sum of irreducibles for $\mathfrak{l}$.

\medskip

This is what we had before! For example, the expression
$
\bigl[\;  \figins{-4}{0.20}{one}\;  \bigr]_{\! D_n}\!\!  =  
\bigl[\;  \figins{-4}{0.20}{one-u}\;  \bigr]_{\! D_n}\!\!  +    
\bigl[\;  \figins{-4}{0.20}{one-d}\;  \bigr]_{\! D_n}
$,
which can be obtained from the extension of Theorem~\ref{thm:jaeg} to tangles,
is a diagrammatic interpretation of the isomorphism  
$V_{fund}(\mathfrak{so}_{2n})\cong V_{fund}(\mathfrak{sl}_{n})\oplus V^*_{fund}(\mathfrak{sl}_{n}) $ 
for $\mathfrak{so}_{2n}\supset\mathfrak{sl}_{n}$, and 
$
\bigl[\  \figins{-4}{0.20}{one}\  \bigr]_{\! B_n}\!\! =   
\bigl[\  \figins{-4}{0.20}{one-u}\  \bigr]_{\! B_n}\!\! +    
\bigl[\  \figins{-4}{0.20}{one-d}\  \bigr]_{\!B_n}\!\! +    
\bigl[\  \figins{-4}{0.20}{one-dash} \ \bigr]_{\! B_n}
$
corresponds to 
$V_{fund}(\mathfrak{so}_{2n+1})\cong V_{fund}(\mathfrak{sl}_{n})\oplus V^*_{fund}(\mathfrak{sl}_{n})\oplus V_{triv}(\mathfrak{sl}_{n})$ 
for $\mathfrak{so}_{2n+1}\supset\mathfrak{sl}_{n}$. 

\medskip

The general picture of \emph{categorification of quantum link invariants}, pioneered by M. Khovanov~\cite{kho} 
lifts the representations $W$ appearing in the RT picture to categories $\mathcal{C}_{\mathfrak{g}}(W)$ 
(which are required to satisfy certain properties) 
and the RT map $f_{RT}$ to a (derived) functor $\mathcal{F}_{RT}$ between (the derived categories of) these categories. 
Again, the isomorphism class of this functor depends only on the isotopy class of the tangle.
The categorification of $f_{RT}$ for ge\-ne\-ral f.d. irreducible representations of QEAs 
was constructed by B. Webster in~\cite{webster1, webster2}. 
\\[1ex]
$$
\mspace{-120mu}
\labellist\small
\pinlabel $W_1^*$    at  38 -30
\pinlabel $\otimes$  at  88 -30
\pinlabel $W_2$     at  152 -30
\pinlabel $\otimes$ at  208 -30
\pinlabel $W_1$     at 270 -30
\pinlabel $W_3^*$    at  25 275
\pinlabel $\otimes$  at  88 275
\pinlabel $W_3$      at 155 275
\pinlabel $\otimes$  at  208 275
\pinlabel $W_2$      at 265 275
\pinlabel $\inv_{\mathfrak{g}}(W,W')$ at 430 75
\pinlabel \begin{turn}{-90}$\in$\end{turn} at 430 135
\pinlabel $f_{RT}$ at 430 185
\endlabellist
\figins{-25}{0.95}{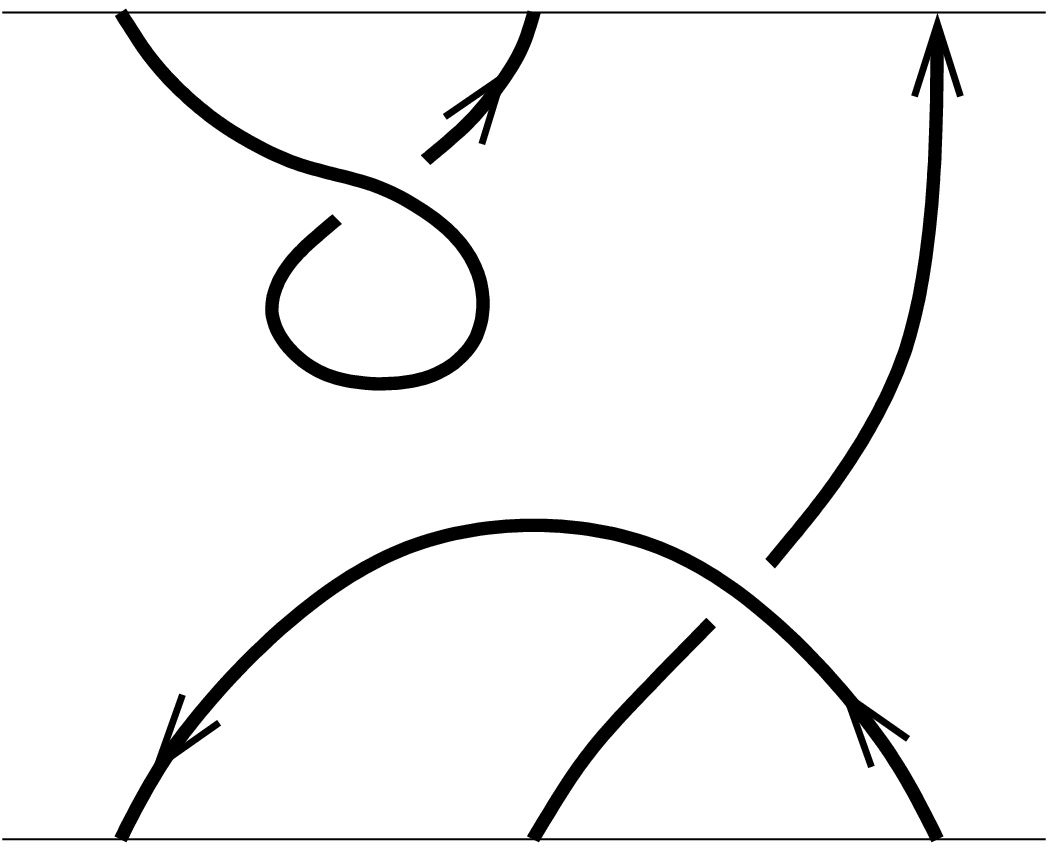}
\mspace{120mu}\text{\scalefont{2}$\leadsto$\normalsize}\mspace{30mu}
\labellist\small
\pinlabel $\mathcal{C}_{\mathfrak{g}}(W_1^*,W_2,W_1)$  at  152  -35
\pinlabel $\mathcal{C}_{\mathfrak{g}}(W_3^*,W_3,W_2)$  at  152  275
\pinlabel $\fun_{\mathfrak{g}}(\mathcal{C}_{\mathfrak{g}}(W),\mathcal{C}_{\mathfrak{g}}(W'))$ at 550 75
\pinlabel \begin{turn}{-90}$\in$\end{turn} at 550 135
\pinlabel $\mathcal{F}_{RT}$ at 550 185
\endlabellist
\figins{-25}{0.95}{tangle2}
\vspace{4ex}
$$

We can try to use Webster's work to 
construct categorical $\mathfrak{l}$-expansions for the categorified $\mathfrak{g}$-RT invariants. 
The categories $\mathcal{C}_{\mathfrak{g}}$ appearing in~\cite{webster1} extend to  
 linear combinations of (arbitrary)  f.d. irreducibles of $\mathfrak{g}$ which  
means that Webster's functors extend to (formal) 
linear combinations of  tangles. 

\begin{defn}
A \emph{categorical Jaeger expansion} consists of (i) \emph{categorified branching rules} i.e. 
a functor 
$
\mathcal{C}_{\mathfrak{g}}(V^{\mathfrak{g}})\to 
\mathcal{C}_{\mathfrak{l}}(\oplus_{i} V_{i}^{\mathfrak{l}}) 
$ for $\mathfrak{l}\subset\mathfrak{g}$, 
which is full, bijective on objects and descends to a map between the respective Grothendieck groups giving   
an isomorphism of $\mathfrak{l}$-representations $V^{\mathfrak{g}}\cong \oplus_{i} V_{i}^{\mathfrak{l}}$,   
and (ii) its extension to corresponding decompositions of the ``tangle functors''. 
Here $V^{\mathfrak{g}}$ and each of the $V_i^{\mathfrak{l}}$
are irreducible f.d. representations of $\mathfrak{g}$ and 
$\mathfrak{l}$ respectively (resp. tensor products of such representations).
\end{defn}

Although (ii) seems desirable from the topological point of view (work still in progress), 
the fulfillment of (i) (see~\cite{vaz1,vaz2}) is already very interesting, 
due to the potential applications to areas like representation theory and physics.

\begin{thm}
There are functors 
$
\mathcal{C}_{\mathfrak{g}}(V^{\mathfrak{g}})\to 
\mathcal{C}_{\mathfrak{l}}(\oplus_{i} V_{i}^{\mathfrak{l}}) 
$ 
categorifying the branching rules for 
$\mathfrak{sl}_{n+1}\supset\mathfrak{sl}_{n}$, 
$\mathfrak{so}_{2n+1}\supset\mathfrak{so}_{2n-1}$,
$\mathfrak{sp}_{2n}\supset\mathfrak{sp}_{2n-2}$,
$\mathfrak{so}_{2n}\supset\mathfrak{so}_{2n-2}$,
(for all finite dimensional representations and tensor products of minuscule representations) ,
$\mathfrak{so}_{2n},\ \mathfrak{so}_{2n+1}\supset\mathfrak{sl}_{n}$ (for fundamental representations and their tensor products).
\end{thm}


\vspace*{1cm}




\begin{thebibliography}{19}


\bibitem{kauffman} 
L.~Kauffman,
\emph{Knots and Physics}, World Scientific, Singapore, 1991.

\bibitem{kho} 
M.~Khovanov,
``A categorification of the Jones polynomial'',
\emph{Duke. Math. J.},  Vol.~101, No.~3 (2000), pp. 359-426.

\bibitem{resh-tur} 
N.~Reshetikhin e V.~Turaev,
``Ribbon graphs and their invariants derived from quantum groups'',
\emph{Commun. Math. Phys.},  Vol.~127, No.~1 (1990), pp. 1-26.

\bibitem{tur} 
V.~Turaev,
``The Yang-Baxter equation and invariants of links'',
\emph{Invent. Math.},  Vol.~92, No.~3 (1988), pp. 527-553.


\bibitem{VW} 
P.~Vaz e E.~Wagner,
``A remark on BMW algebra, $q$-Schur algebras and categorification'', 
\emph{arXiv}:1203.4628v1 [math.QA]  (2012).

\bibitem{VW2} 
P.~Vaz e E.~Wagner,
``(work in progress)'', 
\emph{}  (2012).




\bibitem{webster1} 
B.~Webster,
``Knot invariants and higher representation theory I: diagrammatic and geometric categorification of tensor products'', 
\emph{arXiv}:1001.2020v7 [math.GT]  (2011).


\bibitem{webster2} 
B.~Webster,
``Knot invariants and higher representation theory II: the categorification of quantum knot invariants'', 
\emph{arXiv}:1005.4559v5 [math.GT]  (2011).



\bibitem{vaz1}
P.~Vaz, 
``KLR algebras and the branching rule I: The Gelfand-Tsetlin basis in type $A_n$'',  
\emph{arXiv}:1309.0330v1 [math.RT]  (2013). 
 

\bibitem{vaz2}
P.~Vaz, 
KLR algebras and the branching rule II: The Gelfand-Tsetlin basis in types $B$, $C$ and $D$,\\ 
(2013) (in preparation). 



\end{thebibliography}
\end{document}